\documentclass[oneside]{article}
\usepackage[total={6.6in,8.75in}, top=1.2in, left=0.8in, includefoot]{geometry}
\pdfoutput = 1

\usepackage{setspace}
\usepackage{Style}
\usepackage{pdfpages}

\theoremstyle{theorem}
\newtheorem{lemm}{Lemma}[section]
\newtheorem{theo}[lemm]{Theorem}
\newtheorem{coro}[lemm]{Corollary}

\theoremstyle{definition}
\newtheorem{defi}[lemm]{Definition}
\newtheorem{exam}[lemm]{Example}
\newtheorem{rema}[lemm]{Remark}

\renewcommand{\thefootnote}{\fnsymbol{footnote}}

\title{\Large\scshape\textbf{Necessary conditions for the positivity of Littlewood-Richardson and plethystic coefficients}}
\author{
    Álvaro Gutiérrez \and 
    Mercedes H. Rosas
}

\renewcommand{\thefootnote}{\arabic{footnote}}
\date{23 September 2021}

\usepackage{fancyhdr}

\fancypagestyle{alim}{\fancyhf{}\fancyfoot[R]{\footnotesize
\vspace{-5.67em}
MR was partially supported by
MTM2016-75024-P and FEDER,\\
PID2020-117843GB-I00, and Proyectos I+D+i FEDER Andalucía US-1262169.\\
AG was partially supported by
Junta de Andalucía FQM-333.}}
\begin{document}

\thispagestyle{alim}
\begin{center}\large {\ }
\vspace{1em}
\textbf{\scshape\Large Necessary conditions for the positivity of Littlewood-Richardson and plethystic coefficients}
\\ \vspace{2em}
\renewcommand{\thefootnote}{\fnsymbol{footnote}}
\renewcommand{\arraystretch}{1}
\setlength{\tabcolsep}{2em}
\begin{tabular}{cc}
    \'Alvaro Gutiérrez\footnotemark[1]{}
    & Mercedes H. Rosas\footnotemark[2]{} \\
    {\small Departamento de Álgebra} & {\small Departamento de Álgebra}\\
    {\small Universidad de Sevilla, Spain} & {\small Universidad de Sevilla, Spain}
\end{tabular}
\footnotetext[1]{agutierrez1@us.es}
\footnotetext[2]{mrosas@us.es}
\renewcommand{\thefootnote}{\arabic{footnote}}
\setcounter{footnote}{0}
\\ \vspace{2em}
\end{center}

\begin{abstract}
    We give necessary conditions for the positivity of Littlewood-Richardson coefficients and SXP coefficients. We deduce necessary conditions for the positivity of the plethystic coefficients. Explicitly, our main result states that if $S^\lambda(V)$ appears as a summand in the decomposition into irreducibles of $S^\mu(S^\nu(V))$, then $\nu$'s diagram is contained in $\lambda$'s diagram.
\end{abstract}

{\small
\noindent\hspace{2.75em}\textbf{Keywords:} symmetric functions, plethysm, Littlewood-Richardson coefficients, SXP rule\\

\vspace{-1.2em}
\noindent\hspace{2.75em}\textbf{MSC:} 05E05, 05E18, 05A17}

\tableofcontents
\newpage

\section{Introduction}

The operations of restriction, tensor product, and composition of representations allow us to combine complex representations 
of the general lineal groups, and obtain new  interesting  representations of these groups. 
Breaking these new representations as sums of irreducibles representations is a major problem in representation theory. 
It is in this setting that the families of coefficients that we study in this work appear: the Littlewood--Richardson, Kronecker, and plethystic coefficients, respectively, describe the multiplicities that govern these decompositions. In addition to their importance in representation theory, these coefficients naturally appear in many 
different fields of mathematics from invariant theory, Schubert calculus, and algebraic geometry, to physics and computer science,
\cite{Fulton-Harris, macdonald,  Fulton-Young, MulSo, stanleyEC2}.
For recent work see \cite{COSSZ, Ikenmeyer-Fisher, MiRoSu,pak17, pak20, wildon}.

In the language of symmetric functions, the irreducible representation  of $\mathbf{GL}(V)$ 
indexed by a partition $\lambda$ translates to the Schur function $s_\lambda$.
The tensor product of irreducible representations translates to the ordinary product of Schur functions, 
which allows us to define the Littlewood--Richardson coefficients  as the structural constants for the ordinary product of Schur functions, \(s_\mu\cdot s_\nu = \sum_\lambda c_{\mu, \nu}^\lambda s_\lambda\).
Since Schur polynomials form an orthonormal basis of the space of symmetric functions, we can also write
\(c_{\mu,\nu}^\lambda = \langle s_\lambda, s_\mu\cdot s_\nu\rangle\). 
The operation of composition of representations translates to the plethysm of symmetric functions, which in turn allows us to define the plethystic coefficients \(a_{\mu[\nu]}^\lambda\) 
as the number \(\langle s_\lambda, s_\mu[s_\nu]\rangle\).
Finally, the Kronecker coefficients $g_{\mu,\nu}^\lambda$ are the multiplicities governing the decomposition into irreducibles of the restriction 
of $\mathbf{GL}(V\otimes W)$ to $\mathbf{GL}(V) \times \mathbf{GL}(W)$, via the Kronecker product of matrices.

A famous result, often attributed to Dvir, gives a necessary condition that a Kronecker coefficient must satisfy in order to be nonzero. 
This is a remarkable result, as Dvir's conditions are both elegant and very easy to manipulate. 
Let us identify a partition $\lambda$ with its Ferrers diagram.
Explicitly, given partitions $\mu$ and $\nu$, Dvir defines a rectangular partition $R = (|\mu \cap \nu|^{|\mu \cap \nu'|})$ and shows that if $g_{\mu,\nu}^\lambda$ is nonzero, then $\lambda\subseteq R$ (see \cite{dvir}). 

Dvir's result gives a powerful tool in representation theory. To give just two recent applications,  
Pak and Panova  used it to find a counterexample of the Kirillov-Klyachko conjecture \cite{pak20}, and
Briand, Orellana, and the second author used it to give sharp bounds for the stability of 
the Kronecker products of Schur functions \cite{briand11}.\\

Following the spirit of Dvir's result, we show in Theorem \ref{(c,r)-hooks} that if $s_\lambda$ appears as a nonzero summand on the decomposition of $s_\mu[s_\nu]$ in the Schur basis, then the diagram of $\nu$ is contained in the diagram of $\lambda$. In other words, if $a_{\mu[\nu]}^\lambda$ is nonzero, then $\nu\subseteq\lambda$.

With this aim in mind, we first show in Theorem \ref{Main} how to define a partition $\Theta$, from partitions $\mu$ and $\nu$,  
such that $c_{\mu,\nu}^\lambda\ne0$ implies $\lambda\subseteq\Theta$. Our main tool comes from plethystic calculus, and the operation of evaluation into sums and differences of alphabets. This approach to the study of structural constants has been proven successful in the past for Kronecker \cite{rosas} and plethystic \cite{langley} coefficients.

In Theorems \ref{lowerBound} and \ref{upperBound}, we use this result together with the SXP rule \cite{littlewood, wildon} to determine upper and lower bounds for the partitions $\mu$ appearing with positive coefficient in the expansion of $p_n[s_\lambda]$ in the Schur basis. Explicitly, we show that $\lambda\subseteq\mu\subseteq\Xi$, where $\Xi$ is a purposely crafted partition.
Then, we express the plethysm of two arbitrary Schur functions in terms of SXP coefficients and plethysms of the type $p_n[s_\lambda]$ as follows:
\begin{equation}\label{simplification}
    s_\mu[s_\nu] = \sum_\lambda \frac{\chi^\mu(\lambda)}{z_\lambda} p_\lambda[s_\nu] = \sum_\lambda \frac{\chi^\mu(\lambda)}{z_\lambda} \prod_i p_{\lambda_i}[s_\nu] =
    \sum_\lambda \frac{\chi^\mu(\lambda)}{z_\lambda} \prod_i \left(\sum_\tau b_{\lambda_i[\nu]}^{\tau} s_\tau\right),
\end{equation}
where $\chi^\mu(\lambda)$ is the value of the character $\chi^\mu$ of the Specht module $S^\mu$ on any permutation of cycle type $\lambda$, the number $z_\lambda$ denotes the order of the centralizer of a permutation of cycle type $\lambda$, and $b_{\lambda_i[\nu]}^\tau$ is defined as the number $\langle s_\tau, p_{\lambda_i}[s_\nu]\rangle$.
The tools thus far developed suffice to show our aforementioned main result (Theorem \ref{(c,r)-hooks}).

Finally, in Corollary \ref{trivial,sign}, we completely characterize the multiplicity of the trivial and sign representations on the decomposition into irreducibles of the composition of arbitrary irreducible representations. These last results can also be deduced from Yang’s work \cite{yangFT}.

\section{Preliminaries}\label{sec:preliminaries}

\subsection{Partitions and symmetric functions}

We follow Stanley \cite{stanleyEC2} for the standard concepts and notations in the theory of symmetric functions, the main exception being that we represent our partitions with the French convention\footnote{In the French convention, we use are bottom-left justified diagrams \cite{macdonald}. The coordinate system is cartesian, the origin being aligned with the bottom-left corner.}.

\begin{figure}[h]
    \centering
    \begin{tikzpicture}[x=1em,y=1em, yscale=-1]
    \path[fill=blue!20]
        (0,1)--(1,1)--(1,2)--(2,2)--(2,3)--(4,3)--(4,4)--(0,4);
    \draw[very thick, black]
    (0,0)--(0,4)--(5,4)
    (3,0)--(3,2)--(5,2);
    \draw[very thick,decorate,decoration={calligraphic brace,amplitude=5pt}]
    (0,-2)--(3,-2) node[anchor=south,black,midway,yshift=5pt]{$c$};
    \draw[very thick,decorate,decoration={calligraphic brace,amplitude=5pt}]
    (7,2)--(7,4) node[anchor=west,black,midway,xshift=5pt]{$r$};
    \filldraw[red] (3.5,1.5) node {$\times$};
    \filldraw[black]
    (0,-1.2) node{$\vdots$}
    (3,-1.2) node{$\vdots$}
    (6,4) node{$\ldots$}
    (6,2) node{$\ldots$}
    (1,3) node{\footnotesize$\lambda$}
    (3.5, 1.5) node[anchor=south west]{$(c,r)$};
    \end{tikzpicture}
    \caption{The diagrams of $(c,r)$-hooks fit in the depicted fat-hook region.}
    \label{fig:hook-lookig}
\end{figure}
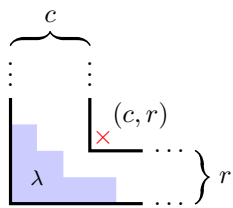

A partition is a weakly decreasing sequence of natural numbers in which there are finitely many nonzero entries. Define the (Ferrers) diagram of a partition $\lambda$ as the subset of $\N_0^2$ made of the points $(c,r)$ such that $0\le c<\lambda_r$. We will often identify a partition with its diagram. A partition is a $(c,r)$-hook if its diagram does not contain the point $(c,r)$ \cite{remmel84}. Note that a $(1,1)$-hook --- usually known just as a hook --- is also a $(c,r)$-hook for any $c\ge1$ and $r\ge1$. We say that $(c,r)$-hooks fit in a fat-hook region of $\N_0^2$ with $c$ columns and $r$ rows (see Figure \ref{fig:hook-lookig}).

\begin{figure}[h]
    \centering
    \begin{tikzpicture}[x=1em,y=1em]
    \path[fill=blue!20]
        (-0.5,-0.5)--(2.5,-0.5)--(2.5,1.5)--(0.5,1.5)--(0.5,2.5)--(-0.5,2.5);
    \ejes{5}{5}
    \foreach \i/\r in   {0/3,1/3,2/1}{
        \foreach \c in {1,...,\r}{
            \filldraw[black] (\c-1,\i) circle (2pt);
        }
    }
    \end{tikzpicture}
    $\quad\quad$
    \begin{tikzpicture}[x=1em,y=1em]
    \path[fill=red!20]
        (5,-0.5)--(2.5,-0.5)--(2.5,1.5)--(0.5,1.5)--(0.5,2.5)--(-0.5,2.5)--(-0.5,5)--(5,5);
    \ejes{5}{5}
    \foreach \p in {(0,3),(1,2),(3,0)}
    {\filldraw[red] \p node {$\times$};}
    \end{tikzpicture}
    \caption{On the left, the Ferrers diagram of $\mu = (3,3,1)$. On the right, its associated ideal is shaded. It is spanned by the set $\OC{\mu}$ of outer corners of $\mu$ (depicted as crosses).}
    \label{fig:OC}
\end{figure}
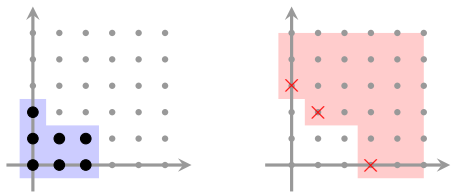

The point $(c,r)$ is an outer corner if $(c,r)\not\in\lambda$ but its addition to the diagram produces a valid partition.
Let $\OC{\lambda}$ be the set of outer corners of $\lambda$. For example, \yindex$\OC{(3,3,1)} = \{(0,3),(1,2),(3,0)\}$. 
The complement $(\lambda)^c$ of the diagram of a partition $\lambda$ defines an ideal
of $\N_0^2$ with respect to the coordinate-wise sum; that is, $(\lambda)^c = \{(x,y)\not\in\lambda\}$ is closed under the sum. The set $\OC{\lambda}$ is the minimal spanning set of $(\lambda)^c=\langle\OC{\lambda}\rangle$. Conversely, the complement of such an ideal containing at least one point of the form $(0,r)$ and one point of the form $(c,0)$ is the diagram of a partition. See Figure \ref{fig:OC} for an illustration of these concepts.

Let the sum of two partitions $\lambda$ and $\mu$ be the partition $\lambda+\mu=(\lambda_1+\mu_1,\lambda_2+\mu_2,...)$, and let the union $\lambda\cup\mu$ of partitions $\lambda$ and $\mu$ be the partition resulting from the sorting of their parts.
Moreover, for a given $n\in\N_+$, we let $n\lambda = \lambda + \lambda + \stackrel{n\text{ times}}{\ldots} + \lambda$ and $\cup^n\lambda = \lambda \cup \lambda \cup \stackrel{n\text{ times}}{\ldots} \cup \lambda$.
We shall write ($\trianglelefteq$) for the dominance order on partitions, letting $\lambda\trianglelefteq\mu$ if $\sum_{i\le j}\lambda_i\le\sum_{i\le j}\mu_i$ for all $j$. We write $\mu\subseteq\lambda$ and say $\mu$ is contained in $\lambda$ whenever $\mu_i \le \lambda_i$ for all $i$. For $\mu\subseteq\lambda$, we let $\lambda/\mu$ denote the set of points in $\lambda$ and not in $\mu$. A rim hook of $\lambda$ is a skew partition $\lambda/\mu$ whose diagram is (orthogonally) connected and contains no $2\times2$ arrangement.
    
Let $n\in\N_+$. The $n$-quotient and $n$-core of a partition encode all of the information of the original partition. We list some of their properties, and refer to Macdonald \cite{macdonald} for their proofs. The $n$-quotient of a partition $\lambda$ is defined as the $n$-tuple $\lambda^* = (\lambda^{(0)},\lambda^{(1)},...,\lambda^{(n-1)}),$
where $\lambda^{(i)}$ is made of the points ($k,j$) in $\lambda$ such that $c_k:=\lambda'_k+k+1\equiv i \ (n)$ and $r_j:=\lambda_{j}+j\equiv i \ (n)$.
Note that $c_k$  only depends on the column and $r_j$ on the row. The $n$-core is defined as the partition $\tilde{\lambda}$ which remains after removing (step by step) every rim hook of length $n$ from $\lambda$. The order in which the rim hooks are removed does not matter. Letting $|\lambda^*|$ be $|\lambda^{(0)}|+|\lambda^{(1)}|+\cdots+|\lambda^{(n-1)}|$, we get the following formula
    \begin{equation}\label{formula sizes}
        |\lambda|=|\tilde{\lambda}|+n|\lambda^*|.
    \end{equation}

Let $\Lambda$ be the algebra of symmetric functions. That is, the algebra $\Q[p_1, p_2, ...]$ spanned by the algebraically independent variables $p_k$ which we name the {power sum symmetric functions}. It will sometimes be useful to identify an element $f$ in $\Lambda$ with a formal power series. Let $X=x_1+x_2+...$ be an alphabet --- a collection of variables called letters. We will identify any $f\in\Lambda$ with its image $f[X]$ under the morphism that maps $p_k$ to $x_1^k + x_2^k +...$ . In particular, we identify $p_1$ with $X$. We write $f[X] = f(x_1, x_2, ...)$ and say that it is the evaluation of $f$ in $X$.

For $f\in\Lambda$, let \(\langle s_\lambda, f\rangle\) denote the coefficient of $s_\lambda$ in the decomposition of $f$ in the Schur basis. Hence \(f=\sum_\lambda \langle s_\lambda, f\rangle \cdot s_\lambda\). We let the set \(\{\lambda : \langle s_\lambda, f\rangle\ne0\}\) be the support of $f$, denoted as \(\supp(f)\). The generalized Littlewood--Richardson coefficient \(c_{\mu^0, \mu^1, ..., \mu^{n-1}}^\lambda\) is defined as the number \(\langle s_\lambda,\ s_{\mu^0}\cdot s_{\mu^1} \cdots s_{\mu^{n-1}}\rangle\). Note that for $n=2$, we recover the usual Littlewood--Richardson coefficient (hereafter, LR coefficient). As an immediate consequence of the Littlewood--Richardson rule, we get the following lemma.
\begin{lemm}\label{LR bounds} If \(\lambda\in\supp(s_\mu\cdot s_\nu)\) then \(\mu\cup\nu\trianglelefteq\lambda\trianglelefteq\mu+\nu\). Moreover, \(c_{\mu,\nu}^{\mu\cup\nu} = c_{\mu,\nu}^{\mu+\nu} = 1\).
\end{lemm}

\subsection{Plethysm}

The notion of plethysm, denoted by $\cdot[\cdot]$, comes from that of composition. Let $f$ and $g$ in $\Lambda$. If $g[X]$ is a sum of monic monomials, $g[X] = g_1 + g_2 + ...$ then $f[g[X]]=f(g_1, g_2, ...)$. In particular, since we identify $p_1$ with $X$, then $f[X]$ is just the plethysm of $f$ with $X$.

\begin{exam}\label{alphabet}
If $f[X]$ is a power series with positive integers as coefficients, it can be expressed as a sum of monic terms. For instance, $2p_2[X] = 2x_1^2 + 2x_2^2 +... = x_1^2 + x_1^2 + x_2^2 + x_2^2 + ...$ .
Consequently,
$p_n[2p_2[X]] = p_n(x_1^2, x_1^2, x_2^2, x_2^2,...) = 2p_{2n}[X].$
\end{exam}


More precisely, the operation of plethysm of symmetric functions is defined axiomatically.

\begin{defi} The plethysm of symmetric functions, denoted by $\cdot[\cdot]$, is the operation $\Lambda\times\Lambda\to\Lambda$ verifying
\begin{enumerate}
    \item $p_n[p_m] = p_{nm}$ for all $n, m\in\N_0$.
    \item For any $f\in\Lambda$, the map $g\mapsto g[f]$ is a $\Z$-algebra homomorphism on $\Lambda$.
    \item For any $f\in\Lambda$, the equality $p_n[f] = f[p_n]$ holds.
\end{enumerate}
\end{defi}

\begin{exam}
We use axiom (3) to compute $p_n[-X] = p_n[-p_1] = -p_1[p_n]$. Then, using axiom (2), we get
\begin{gather*}
s_2[-X] = \frac{p_{1,1}+p_2}{2}[-X] = \frac{1}{2}p_1[-X]p_1[-X] + \frac{1}{2}p_2[X] \\ 
= \frac{1}{2}p_1[X]p_1[X] - \frac{1}{2}p_2[X] = \frac{p_{1,1}-p_2}{2}[X] = s_{1,1}[X].
\end{gather*}
\end{exam}

The core tools of this work come from plethystic calculus. Namely, from the operation of evaluation in sums and differences of alphabets. This next lemma is standard. More general formulas, for $s_\lambda[f\pm g]$ on two arbitrary symmetric functions, are found in \cite{macdonald}.

\begin{lemm}\label{sergeev} Let $X$ and $Y$ be two alphabets and let $\lambda$ be a partition. Then:
\begin{enumerate}
    \item $s_\lambda[-X] = (-1)^{|\lambda|}s_{\lambda'}[X].$
    \item $s_\lambda[X+Y] = \sum_{\mu\subset\lambda} \ s_\mu[X] \ \cdot \ s_{\lambda/\mu}[Y] = \sum_{\mu,\nu}c_{\mu,\nu}^\lambda \ s_\mu[X]\ \cdot \ s_{\nu}[Y].$
    \item\label{sergeev-} $s_\lambda[X-Y] = \sum_{\mu\subset\lambda}(-1)^{|\lambda/\mu|}\ s_\mu[X]\ \cdot \ s_{(\lambda/\mu)'}[Y] = \sum_{\mu,\nu}(-1)^{|\nu|}\ c_{\mu,\nu}^\lambda \ s_\mu[X]\ \cdot \ s_{\nu'}[Y].$
\end{enumerate}
\end{lemm}
\begin{rema}\label{note:positive and negative} Let $X$ and $Y$ be two alphabets. Then, Lemma \ref{sergeev} says that $s_\lambda[X-Y]$ is the generating function of the tableaux $T$ on positive letters from $X$ and negative letters from $-Y$ obeying the semistandarity rules for the positive entries and the opposite rules for the negative ones. 
For instance, in Figure \ref{fig:negative} we have four such tableaux of weights $x_1^2x_2x_3^3x_4$, $(-1)^8y_1y_2y_3y_4^3y_5^2$, $(-1)^3x_1^2x_2x_3^3y_1y_2^2$, and $(-1)^{10}x_1^2x_2^2y_1^3y_2^3y_3^4$, respectively.
\begin{figure}[h]
    \centering
    \ytableaubig
    \footnotesize
    \ytableaushort{4,233,1113}$\quad\quad$
    \ytableaushort{{-1},{-4}{-3},{-5}{-4},{-5}{-4}{-2}}$\quad\quad$
    \ytableaushort{3,{-1},{-2}23,{-2}113}$\quad\quad$
    \ytableaushort{{-3},{-3}{-2}{-1},{-3}{-2}{-1}22,{-3}{-2}{-1}11}
    \ytableausetup{nosmalltableaux}
    \caption{Four valid SSYT with positive and/or negative letters.}
    \label{fig:negative}
\end{figure}
\end{rema}

\begin{nota}\label{note: different evaluations}
In general, evaluating on the alphabet $X+\stackrel{c\text{ times}}{\ldots}+X$ is not equivalent to evaluating on the alphabet $cx_1 + cx_2 + ...$ . We denote the first with $f[cX]$ and the latter with $f[tX]|_{t=c}$. In particular, $-p_k[X] = p_k[-X] \ne \left.p_k[tX]\right|_{t=-1} = (-1)^k p_k[X]$.
\end{nota}

This next theorem enables us to calculate plethysms of the form $p_n[s_\lambda]$.
\begin{theo}[SXP rule \cite{littlewood,wildon}] For any partitions $\lambda, \mu$ and any $n\in\N_+$,
$$\langle s_\mu,\ p_n[s_\lambda]\rangle = \sgn_n(\mu)\cdot \big\langle s_\lambda,\ s_{\mu^{(0)}}\cdot s_{\mu^{(1)}} \cdots  s_{\mu^{(n-1)}}\big\rangle,$$
where $(\mu^{(0)},\ldots,\mu^{(n-1)})$ is the $n$-quotient of $\mu$, and the sign function is defined as in \cite{wildon}.
\end{theo}


\begin{rema} From Equation (\ref{formula sizes}) and the SXP rule, we can deduce that $\mu\in\supp(p_n[s_\lambda])$ implies that $\tilde{\mu}=\emptyset$.
\end{rema}
    
The SXP rule lets us immediately identify some partitions of $\supp(p_n[s_\lambda])$. Let us start with two of them.

\begin{lemm} Let $n\in\N_+$. Then,
\begin{enumerate}[label=(\arabic*), ref={\thelemm.(\arabic*)}]
    \item\label{nlambda} The partition $n\lambda$ is in $\supp(p_n[s_\lambda])$ and $\langle s_{n\lambda}, p_n[s_\lambda]\rangle = 1$.
    \item\label{union lambda} The partition $\cup^n\lambda$ is in $\supp(p_n[s_\lambda])$ and $\langle s_{\cup^n\lambda}, p_n[s_\lambda]\rangle = (-1)^{|\lambda|(n-1)}$.
\end{enumerate}
\end{lemm}
\begin{proof} We prove the first assertion; the second one is shown similarly. Let $\mu=n\lambda$. To begin with, the $n$-core of $\mu$ is empty. Now, checking $\lambda\in\supp(s_{\mu^{(0)}}\ \cdot\ s_{\mu^{(1)}}\ \cdots\  s_{\mu^{(n-1)}})$ will suffice. Since $\lambda=\mu^{(0)}\cup\mu^{(1)}\cup\cdots\cup\mu^{(n-1)}$, the result holds from Lemma \ref{LR bounds}.
\end{proof}

\subsection{A plethystic substitution lemma}
The following lemma links Schur  functions evaluations with LR  coefficients. This result has been used  implicitly  in \cite{langley, rosas}.
Given positive integers $a$ and $b$, let $X_a=x_1+x_2+\ldots+x_a$ and $X_b=x'_1+x'_2+\ldots+x'_b$ be two alphabets. We will identify the alphabet $X_a + X_b$ with $X_{a+b}$ by setting $x'_i\mapsto x_{a+i}$.

\begin{lemm}\label{2.1} Let $\lambda$ be a partition, and let $a,b,c,d \in \N_+$. Then, the evaluation $s_\lambda[X_{a+b}-Y_{c+d}]\ne0$ if and only if there exist partitions $\mu_0$ and $\nu_0$ such that
    $c^{\lambda}_{\mu_0,\nu_0}\ne0$, $\ s_{\mu_0}[X_a-Y_c]\ne0$, and $s_{\nu_0}[X_b-Y_d]\ne0$.
\end{lemm}

\begin{proof}
We know from Lemma \ref{sergeev} that $$s_\lambda[X_{a+b}-Y_{c+d}] = s_\lambda[(X_{a}-Y_{c}) + (X_b -Y_d)] = \sum c_{\mu,\nu}^\lambda\ s_\mu[X_a-Y_c]\  s_\nu[X_b-Y_d].$$
Therefore, if $s_\lambda[X_{a+b}-Y_{c+d}]\ne0$ then at least one of the terms in the sum doesn't vanish.

Conversely, suppose that $c_{\mu_0,\nu_0}^\lambda \ s_{\mu_0}[X_a-Y_c]\ s_{\nu_0}[X_b-Y_d]\ne0$ for some $\mu_0$, $\nu_0$ and consider the equation
$$s_\lambda[X_{a+b}-t Y_{c+d}] = \sum c_{\mu,\nu}^\lambda\ s_\mu[X_a-t Y_c]\ s_\nu[X_b-t Y_d],$$
where $t$ is a variable as in Note \ref{note: different evaluations}. Let $t=-1$. The positivity of LR coefficients ensure that every monomial in both sides of the equality is now positive, so there can't be any cancellation. This means in particular that $s_\lambda[X_{a+b}-t Y_{c+d}]\big|_{t=-1}\ne0$.

On the other hand, we obtain  $s_\lambda[X_{a+b}-t Y_{c+d}]\big|_{t=-1}$ (that we know that is different than zero) from
$s_\lambda[X_{a+b}-Y_{c+d}]$ by setting, for each $j$, the letter $y_j$ to be $-y_j$.
Therefore, we can conclude that $s_\lambda[X_{a+b}-Y_{c+d}]$ is
also different than zero.
\end{proof}

\begin{rema} This result also holds for infinite alphabets. If $\lambda$ is a partition of $n$, working with infinite variables
is equivalent to working with $n$ variables, which is the case
that we settled in the previous lemma.

\end{rema}

\section{Positivity conditions for the Littlewood--Richardson coefficients}\label{sec:LR}
\sectionmark{Positivity conditions for the LR coefficients}
We present a general theorem giving necessary conditions for the positivity of the Littlewood--Richardson coefficients. The following elementary observation will play a crucial role.


\begin{lemm}\label{lem: elementary observation}
We have $s_\lambda[X_r-Y_c]\ne 0$ if and only if $\lambda$ is a $(c,r)$-hook.
\end{lemm}
\begin{proof}
Suppose $s_\lambda[X_r-Y_c]\ne0$. Then, $\lambda$ does not have a $(c,r)$ point, i.e., $(c,r)\in(\lambda)^c$. (In order to see this, think of what would be the value of $(c,r)$ in a tableau with $r$ positive letters and $c$ negative letters.) This, in turn, implies that $\lambda$ fits in a fat-hook region with $r$ rows and $c$ columns (see Figure \ref{fig:hook-lookig}).

Conversely, if $(c,r)$ is in $(\lambda)^c$ then $s_\lambda[X_r-Y_c]\ne0$.
Indeed, the following SSYT in the alphabet $X_r-Y_c$ is always present. Fill the $k$th column with $(k-c)$'s, for $k=0,\ldots,c-1$. After this, the empty cells in the $k$th row are filled with $(k+1)$'s.
\end{proof}

\begin{rema}
In the proof of Lemma \ref{lem: elementary observation}, we showed that if $s_\lambda[X_r - Y_c]\ne 0$ then it fits into a $(c,r)$-hook region. If $\lambda$ fits in multiple of these regions, we can then take the intersection of them to find a smaller region for which $\lambda$ is a subset. See Example \ref{exampleBounds1}.
\end{rema}

\begin{exam}\label{exampleBounds1} Suppose $s_\lambda[X_3-Y_1]\ne0$. This means that there exists a SSYT of shape $\lambda$ and filled with the letters $\{1,2,3,-1\}$, which implies that the point $(1,3)$ does not belong to $\lambda$. Suppose that we also know that $s_\lambda[X_4]\ne0$ and $s_\lambda[-Y_2]\ne0$. This implies that neither $(0,4)$ nor $(2,0)$ belong to $\lambda$. Therefore, $\lambda$ must be a subset of each of the first three regions depicted in Figure \ref{fig:exampleBounds1}. Consequently, $\lambda$ must also be a subset of $(2^3,1)$ (see the fourth diagram in Figure \ref{fig:exampleBounds1}).
\begin{center}
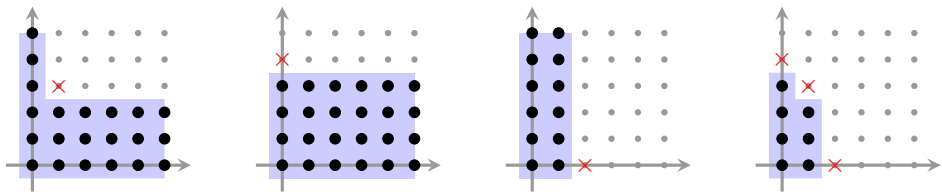
\begin{figure}[h]
    \centering
\begin{tikzpicture}[x=1em,y=1em]
    
    \path[fill=blue!20] (5,-0.5)--(-0.5,-0.5)--(-0.5,5)--(0.5,5)--(0.5,2.5)--(5,2.5);
    \ejes{5}{5}
    \foreach \i/\r in {0/6,1/6,2/6,3/1,4/1,5/1}{
        \foreach \c in {1,...,\r}{
            \filldraw[black] (\c-1,\i) circle (2pt);
        }
    }
    \filldraw[red] (1,3) node {$\times$};
    \end{tikzpicture}$\qquad$
\begin{tikzpicture}[x=1em,y=1em]
    \path[fill=blue!20] (5,-0.5)--(-0.5,-0.5)--(-0.5,3.5)--(5,3.5);
    \ejes{5}{5}
    \foreach \i/\r in {0/6,1/6,2/6,3/6}{
        \foreach \c in {1,...,\r}{
            \filldraw[black] (\c-1,\i) circle (2pt);
        }
    }
    \filldraw[red] (0,4) node {$\times$};
    \end{tikzpicture}$\qquad$
\begin{tikzpicture}[x=1em,y=1em]
    \path[fill=blue!20] (1.5,-.5)--(-0.5,-0.5)--(-0.5,5)--(1.5,5);
    \ejes{5}{5}
    \foreach \i/\r in {0/2,1/2,2/2,3/2,4/2,5/2}{
        \foreach \c in {1,...,\r}{
            \filldraw[black] (\c-1,\i) circle (2pt);
        }
    }
    \filldraw[red] (2,0) node {$\times$};
    \end{tikzpicture}$\qquad$
\begin{tikzpicture}[x=1em,y=1em]
    \path[fill=blue!20]
        (-0.5,-0.5)--(-0.5,3.5)--(0.5,3.5)--(0.5,2.5)--(1.5,2.5)--(1.5,-0.5);
    \ejes{5}{5}
    \foreach \i/\r in {0/2,1/2,2/2,3/1}{
        \foreach \c in {1,...,\r}{
            \filldraw[black] (\c-1,\i) circle (2pt);
        }
    }
    \foreach \p in {(1,3),(0,4),(2,0)}{
    \filldraw[red] \p node {$\times$};}
    \end{tikzpicture}
    \caption{From left to right, the regions $\langle(1,3)\rangle^c$, $\langle(0,4)\rangle^c$, $\langle(2,0)\rangle^c$ and $\langle(1,3), (0,4), (2,0)\rangle^c = (2^3, 1)$.}
    \label{fig:exampleBounds1}
\end{figure}
\end{center}
\end{exam}

In the proof of Lemma \ref{lem: elementary observation}, we constructed a tableau in the alphabet $X_r - Y_c$ for every $(c,r)$-hook. We will refer to it as the canonical SSYT of shape $\lambda$ in the alphabet $X_r-Y_c$ or in the corner $(c,r)$.

\begin{exam}
Let $(c,r) = (3,2)$. Let $\lambda = (5,5,3,1)$. Then, the canonical tableau of shape $\lambda$ in the corner $(c,r)$ is the fourth tableau in Figure \ref{fig:negative}.
\end{exam}

We can now state and prove the main result of this section. For two sets $A$ and $B$ let $A+B = \{a+b\ :\ a\in A, b\in B\}$ be their Minkowski sum.

\begin{theo}\label{Main} Let $n\in\N_+$ and let $\lambda, \mu^0, \mu^1..., \mu^{n-1}$ be partitions. If $c_{\mu^0, \mu^1, ..., \mu^{n-1}}^\lambda$ is nonzero, then \[\lambda\subseteq\Big\langle\sum\limits_{k=0}^{n-1}\OC{\mu^k}\Big\rangle^c,\] where the sum is the Minkowski sum on sets, and the sum in $\N_0^2$ is coordinate-wise.
\comment{If $c_{\mu,\nu}^\lambda\ne0$, then $\overline{R(\mu)} + \overline{R(\nu)} \subseteq \overline{R(\lambda)}$. More generally, if $c_{\mu^0, \mu^1, ..., \mu^{n-1}}^\lambda\ne0$ then $$\sum_k\overline{R(\mu^k)} \subseteq \overline{R(\lambda)}.$$
In other words, if $\lambda\in\supp(s_{\mu^0}\cdot s_{\mu^1}\cdots s_{\mu^{n-1}})$ then $R(\lambda)\subseteq\Theta\left(\sum_k\OC{\mu^k}\right)$.}
\end{theo}

\begin{exam}\label{exam} Let \ytext$\mu^0=(3,2)$, $\mu^1=(1,1) =\mu^2$. Compute their exterior corners:
\begin{center}
    $\mu^0 = \ \ $
    \begin{tikzpicture}[x=0.75em,y=0.75em,baseline=0.5em]
    \foreach \i/\r in {0/3,1/2}{
        \foreach \c in {1,...,\r}{
            \filldraw[black] (\c-1,\i) circle (2pt);
        }
    }
    \foreach \p in {(0,2),(2,1),(3,0)}{
        \filldraw[red] \p node {$\times$};
    }
    \end{tikzpicture}
    then $\OC{\mu^0} = \{(0,2),(2,1),(3,0)\}.$

\vspace{0.5em}
    $\mu^1=\mu^2 = \ \ $
    \begin{tikzpicture}[x=0.75em,y=0.75em,baseline=0.5em]
    \foreach \i/\r in {0/1,1/1}{
        \foreach \c in {1,...,\r}{
            \filldraw[black] (\c-1,\i) circle (2pt);
        }
    }
    \foreach \p in {(0,2),(1,0)}{
        \filldraw[red] \p node {$\times$};
    }
    \end{tikzpicture}
    then $\OC{\mu^1} = \OC{\mu^2} = \{(0,2),(1,0)\}.$
\end{center}
\vspace{0.5em}
Then, add together all possible combinations of exterior corners of our three partitions, to get a set of 9  points of $\N_0^2$,
$\sum_{0}^2\OC{\mu^k} = \{(0,6)$, $(1,4)$, $(2,2)$, $(2,5)$, $(3,3)$, $(4,1)$, $(3,4)$, $(4,2)$, $(5,0)\}.$
The shape $\langle\sum_0^2\OC{\mu^k}\rangle^c = (5, 4, 2, 2, 1, 1)$ arises. Theorem \ref{Main} states that the diagram of every partition in $\supp(s_{\mu^0}s_{\mu^1}s_{\mu^2})$ must be a subset of said region, illustrated in Figure \ref{fig:exam}.
\vspace{1em}
\begin{center}
\begin{figure}[h]
    \centering
\begin{tikzpicture}[x=1em,y=1em]
    \path[fill=blue!20]
        (-0.5,-0.5)--(4.5,-0.5)--(4.5,0.5)--(3.5,0.5)--(3.5,1.5)--(1.5,1.5)--(1.5,3.5)--(0.5,3.5)--(0.5,5.5)--(-0.5,5.5);
    \ejes{6}{7}
    \foreach \i/\r in   {0/5,1/4,2/2,3/2,4/1,5/1}{
        \foreach \c in {1,...,\r}{
            \filldraw[black] (\c-1,\i) circle (2pt);
        }
    }
    \foreach \p in {(0,6),(1,4),(2,2),(2,5),(3,3),(4,1),(3,4),(4,2),(5,0)}
    {\filldraw[red] \p node {$\times$};}
    \end{tikzpicture}
    \caption{The points of $\sum_0^2\OC{\mu^k}$ are depicted as crosses, and the shaded region represents the shape $\langle\sum_0^2\OC{\mu^k}\rangle^c = (5, 4, 2, 2, 1, 1)$.}
    \label{fig:exam}
\end{figure}
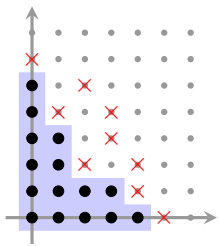
\end{center}
\end{exam}

\begin{proof}[Proof of Theorem \ref{Main}] Let $\lambda\in\supp(s_{\mu^{0}}\cdots s_{\mu^{n-1}})$. Therefore, there exists a partition $\nu^{n-2}$ in $\supp(s_{\mu^{0}}\cdots s_{\mu^{n-2}})$ such that $\lambda\in\supp(s_{\nu^{n-2}}s_{\mu^{n-1}})$. Take now $\nu^{n-2}$. By the same analysis, there exists a partition $\nu^{n-3}$ in $\supp(s_{\mu^{0}}\cdots s_{\mu^{n-3}})$ such that $\nu^{n-2}\in\supp(s_{\nu^{n-3}}s_{\mu^{n-2}})$. Iterate this process to obtain a chain of partitions
$$\lambda=\nu^{n-1}, \nu^{n-2}, ..., \nu^1, \nu^0=\mu^{0}.$$

Choose outer corners $(c_0,r_0)\in\OC{\mu^0}$ and $(c_1,r_1)\in\OC{\mu^1}$. As the canonical tableau for a given corner exists, $s_{\mu^0}[X_{r_0}-Y_{c_0}]\ne0$ and $s_{\mu^1}[X_{r_1}-Y_{c_1}]\ne0$. In addition, we know that $\nu^1\in\supp(s_{\nu^0}s_{\mu^1})=\supp(s_{\mu^0}s_{\mu^1})$. Thus, by Lemma \ref{2.1}, we get $s_{\nu_1}[X_{r_0+r_1}-Y_{c_0+c_1}]\ne0$. Choose now an outer corner $(c_2,r_2)\in\OC{\mu^2}$. Since $\nu^2\in\supp(s_{\nu^1}s_{\mu^2})$, we get that $s_{\nu^2}[X_{r_0+r_1+r_2}-Y_{c_0+c_1+c_2}]\ne0$, again by Lemma \ref{2.1}.

After iterating, $\nu^{n-1}=\lambda$ and so $s_{\lambda}[X_{\Sigma r_i}-Y_{\Sigma c_j}]\ne0$. This means that $\big(\sum c_i, \sum r_j\big)$ is not in $\lambda$. Any choices of corners from $\mu^0, ..., \mu^{n-1}$ will give a similar result, ending the proof.
\end{proof} 

\section{Positivity conditions for the SXP coefficients.}\label{sec:SXP}
\label{Bounds}
In this section we derive necessary conditions for the positivity of the resulting coefficients of the expansion of this plethysm in the Schur basis (SXP coefficients), by combining our previous result for the Littlewood--Richardson coefficients and the SXP rule.

As a corollary of Theorem \ref{LR bounds}, we get the following result.

\begin{theo}\label{lowerBound}
Let $n\in\N_+$ and let $\mu$, $\lambda$ be partitions. If $\langle s_\mu,\ p_n[s_\lambda]\rangle$ is nonzero, then $\lambda\subseteq\mu$.
\end{theo}
\begin{proof} Let $\mu\in\supp(p_n[s_\lambda])$. By the SXP rule, we have $\lambda\in\supp\big(s_{\mu^{(0)}}\cdots s_{\mu^{(n-1)}}\big)$. Choose an outer corner $(c,r)\in\OC{\mu}$. Hence $s_\mu[X_r-Y_c]\ne0$. Let $T:\mu\to\{-c, ..., -2, -1, 1, 2, ..., r\}$
be the canonical SSYT of $\mu$ for $X_r-Y_c$.

Compute the $n$-quotient, thus embedding each $\mu^{(k)}$ inside $\mu$'s diagram. Considering the corresponding values $T(i,j)$ of the canonical tableaux at those embedded cells, we obtain a SSYT $T^k$ of shape $\mu^{(k)}$, which we presume to filled with the alphabet $X_{r_k}-Y_{c_k}$. Then $(c_k,r_k)$ is an outer corner of $\mu^{(k)}$.

Furthermore, we know that no two partitions of the $n$-quotient share any common letters, by construction of $T$ and the $n$-quotient. Consequently, $X_r-Y_c=X_{\Sigma r_k} - Y_{\Sigma c_k}$.

That is, we choose $(c,r)\in\OC{\mu}$ and we show that $(c,r)\in\sum \OC{\mu^{(k)}}$. Therefore, $\OC{\mu}\subseteq\sum \OC{\mu^{(k)}}$. By the Theorem \ref{Main} and the SXP rule, $(c,r)$ is not in $\lambda$.
\end{proof}

\begin{exam}\label{ex:bound1} Let $n=2$, $\lambda=(3,2)$. We have 
\(p_2[s_{3,2}] = -s_{3,3,2,2} + s_{3,3,3,1} + s_{4,2,2,2} - s_{4,3,3} - s_{4,4,1,1} + s_{4,4,2} - s_{5,2,2,1} + s_{5,3,1,1} - s_{5,5} + s_{6,2,2} - s_{6,3,1} + s_{6,4}.\) One can check that 
$(3,2)\subseteq\mu$ for each $\mu\in\supp(p_2[s_{3,2}])$.
\end{exam}

Theorem \ref{lowerBound} gives us a lower bound on the partitions $\mu$ in $\supp(p_n[s_\lambda])$. On the other hand, from the definition of partition we automatically obtain a trivial upper bound. A partition $\mu\in\supp(p_n[s_\lambda])$ must be of size $n|\lambda|$. Hence, the maximum size of the $r$th row is $\big\lfloor\frac{n|\lambda|}{r}\big\rfloor$. 
We refine this upper bound by taking the lower bound into consideration. Let us start with an example.

\begin{exam}\label{ex:bound2}
Let $\mu\in\supp(p_2[s_{3,2}])$. Then, $\mu\subseteq(10,5,3,2,2,1,1,1,1,1)$. We also know $|\mu| = |(2)|\cdot|(3,2)| = 10$. On the other hand, we saw in Example \ref{ex:bound1} that \ytext$(3,2)\subseteq\mu$. However, the only partition of size 10 such that the first row is equal to 10 is the row partition $(10)$, whose diagram clearly does not contain the partition $(3,2)$.
Therefore, our upper bound is subject to improvement.
\end{exam}

\noindent By adjusting our argument, the bound on $\mu_r$ when $\lambda\subseteq\mu$ and $|\mu|= n|\lambda|$ becomes
$$\mu_r\le\left\lfloor\frac{n|\lambda|-|(\lambda_{r+1},\lambda_{r+2},...)|}{r}\right\rfloor=:a_r.$$
A similar analysis for the columns yields the following bounding partition
$$\mu'_c\le\left\lfloor\frac{n|\lambda|-|(\lambda'_{c+1},\lambda'_{c+2},...)|}{c}\right\rfloor=:b_c.$$
Note that these two bounding partitions do not need to be the same. The following result combines both bounding partitions into a more optimized one.
\begin{theo}\label{upperBound} Let $\mu$ and $\lambda$ be two partitions, and let $\Xi^1=(a_1,a_2,...)$ and $\Xi^2=(b_1,b_2,...)'$ with $a_r$ and $b_c$ defined as before. If $\langle s_\mu,\ p_n[s_\lambda]\rangle$ is nonzero, then $\mu\subseteq\Xi^1$ and $\mu\subseteq\Xi^2$. That is, $\mu\subseteq\Xi^1\cap\Xi^2$.
\end{theo}
\begin{exam}\label{ex:bound3} Continuing Examples \ref{ex:bound1} and \ref{ex:bound2}, and by having the lower bound in consideration, we optimize the upper bound to $(8,5,3,2,1,1,1)$. See Figure \ref{fig:bound3}.
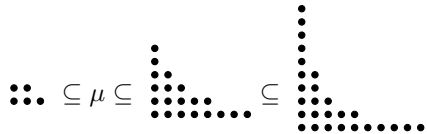
\begin{figure}[h]
    \centering
\begin{tikzpicture}[x=0.5em,y=0.5em,baseline=0em]
    \foreach \i/\r in {0/3,1/2}{
        \foreach \c in {1,...,\r}{
            \filldraw[black] (\c,\i) circle (1.2pt);
        }
    }
    \end{tikzpicture}
$\ \subseteq \mu\subseteq\ $
\begin{tikzpicture}[x=0.5em,y=0.5em,baseline=0.5em]
    \foreach \i/\r in {0/8,1/5,2/3,3/2,4/1,5/1}{
        \foreach \c in {1,...,\r}{
            \filldraw[black] (\c,\i) circle (1.2pt);
        }
    }
    \end{tikzpicture}
$\subseteq\ $
\begin{tikzpicture}[x=0.5em,y=0.5em,baseline=1em]
    \foreach \i/\r in   {0/10,1/5,2/3,3/2,4/2,5/1,6/1,7/1,8/1,9/1}{
        \foreach \c in {1,...,\r}{
            \filldraw[black] (\c,\i) circle (1.2pt);
        }
    }
    \end{tikzpicture}
    \caption{For every $\mu$ in $\supp(p_2[s_{3,2}])$, Examples \ref{ex:bound1}, \ref{ex:bound2} and \ref{ex:bound3} yield the depicted bounds.}
    \label{fig:bound3}
\end{figure}
\end{exam}

\section{Positivity conditions for the general plethystic coefficients}\label{sec:SS}
We now consider the plethysm of two arbitrary Schur functions. Our main result is the following.

\begin{theo}\label{(c,r)-hooks}
Let $\mu$, $\nu$ and $\lambda$ be partitions. If $a_{\mu[\nu]}^\lambda$ is nonzero, then $\nu\subseteq\lambda$.
\end{theo}
\begin{proof}
We shall bring back Equation (\ref{simplification}). 
We know from the LR rule that $\tau, \pi\subseteq\theta$ for all $\theta$ in the support of $s_{\tau}\cdot s_{\pi}$. And if $b_{\lambda_i[\nu]}^\tau$ does not vanish, from Theorem \ref{lowerBound}, we have $\nu\subseteq\tau$. Hence $\nu\subseteq\theta$.
\end{proof}

Two cases of particular interest in representation theory can be further studied. The following corollary can also be deduced from Yang's work \cite{yangFT}. 


\begin{coro}\label{trivial,sign} Let $V$ be a $d$-dimensional vector space, and let  $\mu$ and $\nu$ be partitions of $m$ and $n$ respectively, both of length $\le d$. 
 We have:
\begin{enumerate}[label=(\arabic*), ref={\thecoro.(\arabic*)}]
    \item The coefficient of
    $s_{mn}$ in
    $s_\mu[s_\nu]$ is nonzero if and only if both $\mu$ and $\nu$ are one row partitions. In that case, 
    $\langle s_{mn}, s_{m}[s_{n}]\rangle = 1$.
     \item The coefficient of
    $s_{(1^{mn})}$ in
    $s_\mu[s_\nu]$ is nonzero if and only if both $\mu$ and $\nu$ are one column partitions and $n$ is odd. In that case, 
    $\langle s_{(1^{mn})}, s_{(1^m)}[s_{(1^n)}]\rangle = 1$.
 
\end{enumerate}

\end{coro}
\begin{proof}
Theorem \ref{(c,r)-hooks} shows one implication of each assertion.

From Lemma \ref{nlambda}, 
  $\langle s_{mn}, s_{m}[s_{n}]\rangle =\sum_\lambda \frac{\chi^\mu(\lambda)}{z_\lambda}$, which is the evaluation $s_\mu[1]$. This equals 1 if $\mu$ is a row partition and vanishes otherwise.

On the other hand,  Lemma \ref{union lambda} implies that
 $\langle s_{(1^{mn})}, s_{(1^m)}[s_{(1^n)}]\rangle=\sum_\lambda (-1)^{|\nu|(|\mu|-l(\lambda))}\frac{\chi^\mu(\lambda)}{z_\lambda}$.

If $|\nu|$ is even, this is the evaluation $s_\mu[s_\nu[1]]$, which vanishes unless both $\mu$ and $\nu$ are row partitions. Then $\nu$ must be both a row partition of even size and a column partition. This is impossible.
On the other hand, if $|\nu|$ is odd, the multiplicity of the sign representation is $\sum_\lambda \sgn(\lambda)\frac{\chi^\mu(\lambda)}{z_\lambda} = s_{\mu'}[1]$, which equals 1 if $\mu$ is a column partition and vanishes otherwise.
\end{proof}

The previous lemma can be restated in the language of representation theory as follows.

\begin{coro}
Let $S^\lambda(V)$ be the irreducible representation of $\mathbf{GL}(V)$ indexed by $\lambda$. We have:
\begin{enumerate}[label=(\arabic*), ref={\thecoro.(\arabic*)}]
    \item The trivial representation appears as a summand of $S^\mu(S^\nu(V))$ if and only if $S^\mu$ and $S^\nu$ are trivial representations. In that case, its multiplicity is 1.
    \item The sign representation appears as a summand of $S^\mu(S^\nu(V))$ if and only if both $S^\mu$ and $S^\nu$ are sign representations and $|\nu|$ is odd. In that case, its multiplicity is 1.
\end{enumerate}
\end{coro}

\section{Final Remarks}\label{sec:final}

Recently, there has been plenty of interest in the closely related problem of  understanding  the complexity of deciding whether Littlewood--Richardson \cite{burgisser}, Kronecker \cite{briand09, pak17}, and plethystic \cite{ikenmeyer} coefficients are nonzero. Moreover, and in the cases where it is possible, polynomial algorithms have been developed to determine the positivity of such coefficients.

In \cite{chen}, the plethysm $s_{1,1}[s_{4,2,2}]$ is used as and example to illustrate the importance of this problem, in the case of the plethystic coefficients. A priori, there are $p(16) = 231$ partitions that could appear in $\supp(s_{1,1}[s_{4,2,2}])$, but only 40 actually do. By Theorem \ref{(c,r)-hooks}, we bring this initial number to 142 making it more approachable from the computational perspective.

\begin{figure}[h]
    \centering
    \begin{tikzpicture}[x=1em, y=1em]
    \draw[black!40, arrows=-{\tip},very thick](-1,0)--(17,0);
    \draw[black!40, arrows=-{\tip},very thick](0,-1)--(0,6);
    \foreach \i in {5, 10, 15}{
    \draw (\i,0) -- (\i,-2pt) node[anchor = north]{$\i$};}
    \draw (18, 0) node{$n$};
    \foreach \i in {0.25, 0.50, 0.75, 1.00}{
    \draw (0,\i*5) -- (-2pt,\i*5) node[anchor = east]{$\i$};}
    \foreach \i/\j in {1/1, 2/0.90909090909090, 3/0.833333333333333, 4/0.753246753246753, 5/0.693181818181818, 6/0.633766233766234, 7/0.582070707070707, 8/0.534603174603175, 9/0.492691029900332, 10/0.453961456102784, 11/0.419501133786848, 12/0.387884519107749, 13/0.359371492704826, 14/0.333283183510738, 15/0.309601274485606}{
    \filldraw (\i,5*\j) circle (2pt);}
    \foreach \i/\j in {1/0.333333333333333, 2/0.272727272727273, 3/0.233333333333333, 4/0.155844155844156, 5/0.107954545454545, 6/0.0701298701298701, 7/0.0467171717171717, 8/0.0304761904761905}{
    \draw (\i,5*\j) circle (2pt);}
    \end{tikzpicture}
    \caption{In black, the ratio of partitions of $3n$ that are of length $\ge n$; in white, the ratio of partitions of $3n$ that are in $\supp(s_{2,1}[s_{(1^n)}])$. The high computational cost of plethysm only allowed us to gather data for $n\le 8$.}
    \label{fig: final}
\end{figure}
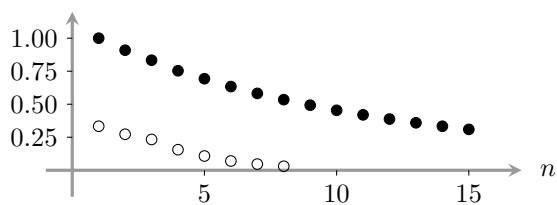

We provide one further example. Let $\mu = (2,1)$, let $\nu = (1^n)$. Any partition in the support of $s_\mu[s_\nu]$ must be, by our theorem, a partition of $3n$ with length at least $n$. This turns out to be a fairly restrictive condition for large $n$. See Figure \ref{fig: final}.

In addition, the theoretical insights of our main theorem are interesting on their own. We hope that the elegant nature of our results will serve as a useful lemmas, and contribute to these complexity results.

\section*{Acknowledgements}

The authors would like to express their gratitude towards Adrià Lillo, Emmanuel Briand, and Laura Colmenarejo for their useful remarks and interesting conversations.

\bibliographystyle{CRBibStyle}
\bibliography{BibliographyPositivity} 

\def\bysame{\leavevmode ---------\thinspace}
\makeatletter
\gdef\og{``}\gdef\fg{''}\makeatother
\def\cdrandname{\&}
\providecommand\cdrnumero{no.~}
\providecommand{\cdredsname}{eds.}
\providecommand{\cdredname}{ed.}
\providecommand{\cdrchapname}{chap.}
\providecommand{\cdrmastersthesisname}{Memoir}
\providecommand{\cdrphdthesisname}{PhD Thesis}
\begin{thebibliography}{10}

\bibitem{briand09}
E.~Briand, R.~Orellana, M.~Rosas, {\og Reduced {K}ronecker coefficients and
  counter-examples to {M}ulmuley's strong saturation conjecture {SH}\fg},
  \emph{Comput. Complexity} \textbf{18} (2009), \cdrnumero 4, p.~577-600,
  \url{https://doi.org/10.1007/s00037-009-0279-z}, With an appendix by Ketan
  Mulmuley.

\bibitem{briand11}
\bysame , {\og The stability of the {K}ronecker product of {S}chur
  functions\fg}, \emph{J. Algebra} \textbf{331} (2011), p.~11-27,
  \url{https://doi.org/10.1016/j.jalgebra.2010.12.026}.

\bibitem{burgisser}
P.~B\"{u}rgisser, C.~Ikenmeyer, {\og Deciding positivity of
  {L}ittlewood-{R}ichardson coefficients\fg}, \emph{SIAM J. Discrete Math.}
  \textbf{27} (2013), \cdrnumero 4, p.~1639-1681,
  \url{https://doi.org/10.1137/120892532}.

\bibitem{chen}
Y.~M. Chen, A.~M. Garsia, J.~Remmel, {\og Algorithms for plethysm\fg}, Contemp.
  Math., vol.~34, 1984, p.~109-153,
  \url{https://doi.org/10.1090/conm/034/777698}.

\bibitem{COSSZ}
L.~Colmenarejo, R.~Orellana, F.~Saliola, A.~Schilling, M.~Zabrocki, {\og The
  mystery of plethysm coefficients\fg}, 2022,
  \url{https://arxiv.org/abs/2208.07258}.

\bibitem{dvir}
Y.~Dvir, {\og On the {K}ronecker product of {$S_n$} characters\fg}, \emph{J.
  Algebra} \textbf{154} (1993), \cdrnumero 1, p.~125-140,
  \url{https://doi.org/10.1006/jabr.1993.1008}.

\bibitem{Ikenmeyer-Fisher}
N.~Fischer, C.~Ikenmeyer, {\og The computational complexity of plethysm
  coefficients\fg}, \emph{Comput. Complexity} \textbf{29} (2020), \cdrnumero 2,
  p.~Paper No. 8, 43, \url{https://doi.org/10.1007/s00037-020-00198-4}.

\bibitem{Fulton-Young}
W.~Fulton, \emph{Young tableaux}, London Mathematical Society Student Texts,
  vol.~35, Cambridge University Press, Cambridge, 1997, With applications to
  representation theory and geometry, x+260~pages.

\bibitem{Fulton-Harris}
W.~Fulton, J.~Harris, \emph{Representation theory}, Graduate Texts in
  Mathematics, vol. 129, Springer-Verlag, New York, 1991,
  \url{https://doi.org/10.1007/978-1-4612-0979-9}, A first course, Readings in
  Mathematics, xvi+551~pages.

\bibitem{ikenmeyer}
C.~Ikenmeyer, K.~D. Mulmuley, M.~Walter, {\og On vanishing of {K}ronecker
  coefficients\fg}, \emph{Comput. Complexity} \textbf{26} (2017), \cdrnumero 4,
  p.~949-992, \url{https://doi.org/10.1007/s00037-017-0158-y}.

\bibitem{langley}
T.~M. Langley, J.~B. Remmel, {\og The plethysm {$s_\lambda[s_\mu]$} at hook and
  near-hook shapes\fg}, \emph{Electron. J. Combin.} \textbf{11} (2004),
  \cdrnumero 1, p.~Research Paper 11, 26,
  \url{http://www.combinatorics.org/Volume_11/Abstracts/v11i1r11.html}.

\bibitem{littlewood}
D.~E. Littlewood, {\og Modular representations of symmetric groups\fg},
  \emph{Proc. Roy. Soc. London Ser. A} \textbf{209} (1951), p.~333-353,
  \url{https://doi.org/10.1098/rspa.1951.0208}.

\bibitem{macdonald}
I.~G. Macdonald, \emph{Symmetric functions and {H}all polynomials}, second
  \cdredname, Oxford Classic Texts in the Physical Sciences, 2015, With
  contribution by A. V. Zelevinsky and a foreword by Richard Stanley, Reprint
  of the 2008 paperback edition [ MR1354144], xii+475~pages.

\bibitem{MiRoSu}
M.~Mishna, M.~Rosas, S.~Sundaram, {\og Vector partition functions and
  {K}ronecker coefficients\fg}, \emph{J. Phys. A} \textbf{54} (2021),
  \cdrnumero 20, p.~Paper No. 205204, 29,
  \url{https://doi.org/10.1088/1751-8121/abf45b}.

\bibitem{MulSo}
K.~D. Mulmuley, M.~Sohoni, {\og Geometric complexity theory. {I}. {A}n approach
  to the {P} vs. {NP} and related problems\fg}, \emph{SIAM J. Comput.}
  \textbf{31} (2001), \cdrnumero 2, p.~496-526,
  \url{https://doi.org/10.1137/S009753970038715X}.

\bibitem{pak17}
I.~Pak, G.~Panova, {\og On the complexity of computing {K}ronecker
  coefficients\fg}, \emph{Comput. Complexity} \textbf{26} (2017), \cdrnumero 1,
  p.~1-36, \url{https://doi.org/10.1007/s00037-015-0109-4}.

\bibitem{pak20}
\bysame , {\og Breaking down the reduced {K}ronecker coefficients\fg}, \emph{C.
  R. Math. Acad. Sci. Paris} \textbf{358} (2020), \cdrnumero 4, p.~463-468,
  \url{https://doi.org/10.5802/crmath.60}.

\bibitem{remmel84}
J.~B. Remmel, {\og The combinatorics of {$(k,l)$}-hook {S}chur functions\fg},
  Contemp. Math., vol.~34, 1984, p.~253-287,
  \url{https://doi.org/10.1090/conm/034/777704}.

\bibitem{rosas}
M.~Rosas, {\og The {K}ronecker product of {S}chur functions indexed by two-row
  shapes or hook shapes\fg}, \emph{J. Algebraic Combin.} \textbf{14} (2001),
  \cdrnumero 2, p.~153-173, \url{https://doi.org/10.1023/A:1011942029902}.

\bibitem{stanleyEC2}
R.~P. Stanley, \emph{Enumerative combinatorics. {V}ol. 2}, Cambridge Studies in
  Advanced Mathematics, vol.~62, 1999,
  \url{https://doi.org/10.1017/CBO9780511609589}, With a foreword by Gian-Carlo
  Rota and appendix 1 by Sergey Fomin, xii+581~pages.

\bibitem{wildon}
M.~Wildon, {\og A generalized {SXP} rule proved by bijections and
  involutions\fg}, \emph{Ann. Comb.} \textbf{22} (2018), \cdrnumero 4,
  p.~885-905, \url{https://doi.org/10.1007/s00026-018-0409-x}.

\bibitem{yangFT}
M.~Yang, {\og The first term in the expansion of plethysm of {S}chur
  functions\fg}, vol. 246, 2002, Formal power series and algebraic
  combinatorics (Barcelona, 1999), p.~331-341,
  \url{https://doi.org/10.1016/S0012-365X(01)00266-7}.

\end{thebibliography}

\end{document}